\documentclass[letterpaper,12pt]{amsart}
\title{Combinatorial Symbolic Powers }

\author{ Seth Sullivant}
\address{Department of Mathematics and Society of Fellows \\ Harvard University, Cambridge, MA 02138}
\email{seths@math.harvard.edu }

\date{}

\usepackage{latexsym,array,delarray,amsthm,amssymb,epsfig}

\theoremstyle{plain}
\newtheorem{thm}{Theorem}[section]
\newtheorem{lemma}[thm]{Lemma}
\newtheorem{prop}[thm]{Proposition}
\newtheorem{cor}[thm]{Corollary}

\theoremstyle{definition}
\newtheorem{defn}[thm]{Definition}
\newtheorem{ex}[thm]{Example}

\newtheorem{ques}[thm]{Question}

\theoremstyle{remark}
\newtheorem*{rmk}{Remark}

\newcommand{\nn}{\mathbb{N}}
\newcommand{\pp}{\mathbb{P}}

\newcommand{\kk}{\mathbb{K}}

\newcommand{\bfa}{\mathbf{a}}
\newcommand{\bfb}{\mathbf{b}}
\newcommand{\bfc}{\mathbf{c}}
\newcommand{\bfd}{\mathbf{d}}
\newcommand{\bfe}{\mathbf{e}}

\newcommand{\bfn}{\mathbf{n}}
\newcommand{\bfp}{\mathbf{p}}
\newcommand{\bfr}{\mathbf{r}}

\newcommand{\bfx}{\mathbf{x}}
\newcommand{\bfy}{\mathbf{y}}
\newcommand{\bfz}{\mathbf{z}}

\newcommand{\ca}{\mathcal{A}}

\newcommand{\cg}{\mathcal{G}}

\newcommand{\ind}{\mbox{$\perp \kern-5.5pt \perp$}}

\begin{document}
\maketitle

\begin{abstract}
Symbolic powers are studied in the combinatorial context of monomial ideals.  When the ideals are generated by quadratic squarefree monomials, the generators of the symbolic powers are obstructions to vertex covering in the associated graph and its blowups.  As a result, perfect graphs play an important role in the theory, dual to the role played by perfect graphs in the theory of secants of monomial ideals.  We use Gr\"obner degenerations as a tool to reduce questions about symbolic powers of arbitrary ideals to the monomial case.  Among the applications are a new, unified approach to the Gr\"obner bases of symbolic powers of determinantal and Pfaffian ideals.

\end{abstract}


\section{Introduction}

The $r$-th symbolic power of an ideal $I$ in a N\"otherian ring $R$ is the ideal
$$I^{(r)} \quad = \quad \left( R^{-1}_I \cdot  I^r \right) \cap R,$$
where $R_I$ denotes the complement of the minimal primes of $I$.
In the down-to-earth setting where $I$ is a radical ideal in a polynomial ring $\kk[\bfx] = \kk[x_1, \ldots, x_n]$ over an algebraically closed field, Zariski and Nagata showed that this is the same operation as the differential power of $I$:
$$I^{<r>} \quad = \quad  \bigcap_{ \bfp \in V(I) }   m_\bfp^r$$
where the intersection runs over all maximal ideals $m_\bfp$ containing $I$ (see e.g.~\cite{Eisenbud1979}).
\begin{thm}[Nagata, Zariski]
If $I$ is a radical ideal in a polynomial ring over an algebraically closed field then
$$I^{(r)} = I^{< r >}. $$
\end{thm}
In characteristic zero, the differential power can also be computed by taking derivatives:
$$I^{< r >} \quad  = \quad  \left<  f   \, \, |   \, \,  \frac{\partial^{|\bfa|}f}{\partial x^{\bfa}}   \in I  \mbox{ for all } \bfa \in \nn^n \mbox{ with } |\bfa| =  \sum_{i = 1}^n a_i \leq r-1  \right>. $$
Thus, the symbolic power $I^{(r)}$ contains all polynomials that vanish to order $r$ on the affine variety $V(I)$, and hence contains important geometric information about the variety.  Among the other applications of symbolic powers are their connections to secant varieties, which was the original motivation for this work.

Our goal in this paper is to study the symbolic powers $I^{(r)}$ for combinatorially defined ideals, and in particular, for squarefree monomial ideals.  One reason for focusing on the monomial case is that we can often bootstrap computations of symbolic powers of monomial ideals to other combinatorially defined ideals.  In particular, we use Gr\"obner degenerations as a tool to reduce questions about symbolic powers of arbitrary ideals to symbolic powers of initial ideals.  This strategy is particularly successful in the case when $I$ is a determinantal or Pfaffian ideal, and provides a new framework for proving many of the classical results about symbolic powers of such ideals (e.g. in \cite{Baetica1998, Bruns2001, Bruns2003, Conca1998}).  This paper should be read as a companion paper to \cite{Sturmfels2006}, extending and exploiting the strategy described there from secant ideals to symbolic powers.

The outline of this paper is as follows.  In the next section, we describe some preliminary results and definitions regarding symbolic powers, and their relations to secant ideals.  We define differentially perfect ideals, which are those ideals whose symbolic powers satisfy a natural recurrence relation.  Sections 3 and 4 are concerned with studying generators for the symbolic powers of edge ideals and antichain ideals, two classes of squarefree monomial ideals of special significance in combinatorial commutative algebra.  

In Sections 5 and 6 we show how the results of Sections 3 and 4 concerning the monomial case can be exploited to prove theorems about the symbolic powers of combinatorially defined ideals, using Gr\"obner degenerations.  Section 5 concerns classical determinantal ideals (of generic and symmetric matrices) and Pfaffian ideals, exploiting some results from \cite{Sturmfels2006}.  Section 6 is concerned with more detailed proofs for some special examples of Segre-Veronese varieties.  We give a new proof of some Gr\"obner basis results for minors of Hankel matrices, and provide two new examples of classes of determinantal ideals whose secants and symbolic powers are well-behaved.  The second of these examples is significant, because the relevant initial ideals are not antichain ideals, and provide examples that do not appear to be amenable to the use of the Knuth-Robinson-Schensted (KRS) correspondence.

To close the Introduction, we give an example to illustrate how the symbolic powers of initial ideals can be used as a  tool to deduce the equations and Gr\"obner bases of the symbolic powers of classical ideals.  Let $V$ denote the Segre embedding of $\pp^1 \times \pp^1 \times \pp^1$ in $\pp^7$.  The ideal $I = I(V)$ is generated by nine quadrics 
$$\underline{x_{001}x_{110}}-x_{100}x_{011}, \quad \underline{x_{010}x_{101}}-x_{100}x_{011}, \quad \underline{x_{111}x_{100}}-x_{101}x_{110}, $$
$$\underline{x_{111}x_{010}}-x_{011}x_{110}, \quad \underline{x_{111}x_{001}}-x_{011}x_{101}, \quad \underline{x_{000}x_{110}}-x_{010}x_{100}, $$
$$\underline{x_{000}x_{101}}-x_{001}x_{100}, \quad \underline{x_{000}x_{011}}-x_{001}x_{010}, \quad  \underline{x_{000}x_{111}}-x_{100}x_{011}$$
that form a Gr\"obner basis with respect to the lexicographic term order $\prec$ with
$$x_{000} \succ x_{111} \succ x_{001} \succ x_{010} \succ x_{100} \succ x_{011} \succ x_{101} \succ x_{110}$$
where the underlined terms are the leading terms.  The initial ideal is the edge ideal $I(G)$ of the graph $G$ with eight vertices and nine edges given by the nine underlined terms of the given binomials.  This graph is bipartite and thus the secant ideal $I(G)^{\{r\}} = \left< 0 \right>$ for $r > 1$.  This implies that the term order $\prec$ is delightful, as defined in \cite{Sturmfels2006}.  Since bipartite graphs are perfect, we deduce by Corollary \ref{cor:delightperfect} that the symbolic powers of $I$ equal the ordinary powers: $I^{(r)}  =  I^r$ for all $r$.  Furthermore, the set of all products of $r$ of the nine quadrics above form a Gr\"obner basis for the symbolic powers $I^{(r)}$ with respect to the given lexicographic term order.
\medskip

\noindent {\bf Acknowledgments.}
I am grateful to Jessica Sidman, whose project to understand the algebraic underpinnings of prolongations \cite{Sidman2006}, and the resulting discussions, led me to the study of symbolic powers.
I also thank Aldo Conca and Rafael Villarreal for useful comments on an earlier version of the paper.


\section{Preliminaries}

In this section, we outline some of the preliminary statements we will need about symbolic powers.  In particular, we develop the relationship between symbolic powers and secant ideals.  As we will often need to exploit the equivalence between symbolic powers and differential powers, we will assume throughout that $\kk$ is an algebraically closed field.  One of the main definitions in this section is the definition of a differentially perfect ideal.  We also give a formula for computing symbolic powers of arbitrary radical ideals in terms of joins.

If $I$ and $J  \subset \kk[\bfx]$ are two ideals, their join is the new ideal
$$I\ast J \quad = \quad    \left( I(\bfy) + J(\bfz) +  \left< x_i - y_i - z_i \, \, | \, \, i = 1, \ldots, n \right>  \right)  \bigcap \kk[\bfx]$$
where $I(\bfy)$ is the ideal $I$ with variable $y_i$ substituted for $x_i$.  The secant ideal $I^{\{2\}}$ is the join of $I$ with  itself:  $I^{\{2\}} = I \ast I$.  The $r$th secant ideal $I^{\{r\}}$ is the $r$-fold join
$$I^{\{r\}} \quad = \quad  I \ast I \ast \cdots \ast I.$$
If $I$ and $J$ are homogeneous radical ideals with varieties $V = V(I)$ and $W = V(J)$, the ideal $I \ast J$ is vanishing ideal of the embedded join
$$V \ast W  \quad = \quad    \overline{   \cup_{v \in V}  \cup_{w \in W}    \left< v,w \right>  }$$
where $\left< v, w \right>$ is the line spanned by $V$ and $W$ and the closure operation is the Zariski closure.  One of the first results relating symbolic powers and secant ideals is a theorem of Catalano-Johnson \cite{Catalano2001}.  Let $m = \left<x_1, \ldots, x_n \right>$ be the homogeneous maximal ideal.

\begin{thm}
Suppose that ${\rm char} \, \kk = 0$ and that $I \subseteq m^2$ is a homogeneous radical ideal.  Then
$$I^{\{r\}}  \subseteq  I^{(r)}.$$
\end{thm}  

By the end of this section, we will provide a proof of the following more general result, which holds over an arbitrary algebraically closed field.

\begin{prop}\label{prop:secdiff}
Let $I$ be a homogeneous radical ideal such that $I \subseteq m^2$.  Then
$$I^{\{r+s -1\}}  \subseteq  \left(I^{\{r\}} \right)^{(s)}.$$
\end{prop}

Besides merely the containment between the secant ideals and symbolic powers, it is known that some graded pieces of the secant ideals and symbolic powers are the same, and thus that, in characteristic zero, some information about the secant ideals $I^{\{r\}}$ can be determined by computing derivatives.

\begin{prop}[\cite{Landsberg2003}, \cite{Sidman2006}]\label{prop:lands}
Suppose ${\rm char} \, \kk = 0$.
Let $I$ be a homogeneous radical ideal such that ${\rm indeg}(I) = k$.  Then the $r(k-1) +1$ graded piece of the $(r-1)(k-1) +1$ symbolic power of $I$ equals the $r(k-1) + 1$ graded piece of the $r$-th secant of $I$:
$$I^{((r-1)(k-1)+1)}_{r(k-1)+1} = I^{\{r\}}_{r(k-1)+1}.$$
\end{prop}

Here ${\rm indeg}(I)$ is the initial degree of $I$, which is the smallest degree of a nonzero polynomial in $I$.  In \cite{Landsberg2003} Proposition \ref{prop:lands} is stated in terms of prolongations in the special case where ${\rm indeg}(I) = 2$.  The prolongation is merely a differential geometry operation identical to taking a particular graded piece of the symbolic power, as shown in \cite{Sidman2006}. 

Of course, it is not possible that $I^{(r)} = I^{\{r\}}$ since they define different varieties.  In particular, for any $i$  such that $1 \leq i \leq r-1$, we always have $I^{(i)}I^{(r-i)} \subset I^{(r)}$.  Thus, for all $r$ and homogeneous $I$ not containing linear forms, 
\begin{equation} 
I^{\{r\}} +  \sum_{i = 1}^{r-1}  I^{(i)} I^{(r-i)} \subseteq  I^{(r)}.
\end{equation} 
More generally, we have the containment:
\begin{equation} \label{eq:diffperfect}  
I^{\{r+s-1\}} +  \sum_{i = 1}^{s-1}  \left(I^{\{r\}}\right)^{(i)} \left( I^{\{r\}} \right)^{(s-i)} \subseteq  \left(I^{\{r\}}\right)^{(s)}.
\end{equation}
For many interesting families of ideals, the containment in (\ref{eq:diffperfect}) is an equality.  This suggests to us the following definition.

\begin{defn}
An ideal $I$ is $r$-\emph{differentially perfect} if for all $s$
$$\left(I^{\{r\}}\right)^{(s)}  = I^{\{r+s-1\}} +  \sum_{i = 1}^{s-1}  \left(I^{\{r\}}\right)^{(i)} \left( I^{\{r\}} \right)^{(s-i)}.$$
An ideal is \emph{differentially perfect} if it is $r$-differentially perfect for all $r$.
\end{defn}

Note that an equivalent definition of $r$-differentially perfect is that the symbolic powers of the secant ideal $I^{\{r\}}$ satisfy:
$$\left(I^{\{r\}}\right)^{(s)} =  \sum_{\lambda \vdash s}   I^{\{r + \lambda_1 -1\}}   I^{\{r + \lambda_2 -1\}} \cdots  I^{\{r + \lambda_{l(\lambda)} -1 \}}  $$
where the sum runs over all partitions $\lambda$ of $s$, and $l(\lambda)$ is the number of parts of $\lambda$.

The simplest ideals with respect to computing symbolic powers are the ones that satisfy $I^{(r)} = I^r$ for all $r$.   Such  ideals are called \emph{normally torsion free} because of their connections to the Rees algebra.  The normally torsion free squarefree monomial ideals were classified by a remarkable result of Gitler, Reyes, and Villarreal \cite{Gitler2006}.  They showed that a squarefree monomial ideal is normally torsion free if and only if the corresponding hypergraph satisfies the max-flow min-cut property.  
Their result makes a strong connection between commutative algebra and combinatorial optimization.

The differentially perfect ideals represent the next simplest possible ideals with respect to computing symbolic powers.  One goal of this paper is to provide the beginnings of a possible classification of differentially perfect ideals.  In particular, in Section \ref{sec:graph}, we classify the $1$-differentially perfect ideals generated by quadratic monomials.  Note that differentially perfect differs from the usual notion of perfect in commutative algebra, which means that the quotient $R/I$ is Cohen-Macaulay.

To provide the proof of Proposition \ref{prop:secdiff}, we need a number of auxiliary results about  joins and secant ideals and the relations to symbolic powers.

\begin{lemma}\label{lem:degreebound}
Let $I$ be a homogeneous radical ideal such that $I \subseteq m^2$.  Then $I^{\{r\}} \subseteq m^{r+1}$.
\end{lemma}

\begin{proof}
We proceed by induction on $r$.
Suppose that $f \in I^{\{r\}}$ is a polynomial of degree less than $r+1$.  Since $I^{\{r\}}  \subseteq I^{\{r-1\}}$, we can assume that $f$ has degree $r$.  Furthermore, since $I^{\{r\}}$ is also radical (\cite{Simis2000}) we can suppose that $f$ is not a power of another polynomial.  Since $f \in I^{\{r\}}$, we have that $f(\bfy + \bfz)  \in I^{\{r-1\}}(\bfy) + I(\bfz)$.  Let $\bfx^\bfa$ be a monomial appearing in $f$.  Plugging in $\bfy +\bfz$, we get the expansion
$$(\bfy + \bfz)^\bfa  =  \sum_{ \bfb + \bfc = \bfa }  \prod_{i =1}^n   {a_i \choose b_i}   y_i^{b_i} z_i^{c_i}.$$
Paying particular attention to the case where $\bfc$ is a standard unit vector $\bfe_i$, we produce the term
$a_i \bfy^{\bfb - \bfe_i}  z_i$.  Since $\bfy^{\bfb-\bfe_i}$ has degree $r-1$ and by our induction hypothesis, this monomial could not appear in any polynomial in  $I^{\{r-1\}}(\bfy) + I(\bfz)$. Thus, the coefficient of this monomial in $f$ is zero, or ${\rm char} \,  \kk$ divides $a_i$ for all $i$.  If $\kk$ has characteristic zero we are done, so suppose that $\kk$ has positive characteristic $p$.  Now we know that $f \in \kk[x_1^p, \ldots, x_n^p]$.  However, we have assumed that $\kk$ is algebraically closed, so there is a polynomial $g \in \kk[\bfx]$ such that $g^p = f$.  This contradicts our assumption that $f$ was not a power.  
\end{proof}

\begin{lemma}\label{lem:intersect}
The join distributes over intersections:
$$ \left( \bigcap_{l \in \mathcal{L}}  J_l  \right)  \ast K  \quad = \quad  \bigcap_{l \in \mathcal{L}} ( J_l \ast K).$$
\end{lemma}

\begin{proof}
A polynomial $f$ belongs to $( \cap  J_l  )  \ast K$ if and only if $f( \bfy + \bfz) \in (\cap J_l)(\bfy) + K(\bfz)$ if and only if $f(\bfy + \bfz) \in J_l(\bfy) + K(\bfz)$ for all $l \in \mathcal{L}$ if and only if $f \in \cap (J_l \ast K)$.
\end{proof}

\begin{lemma}\cite{Simis2000,Sturmfels2006}\label{lem:injoin}
For any term order $\prec$ the join of initial ideals contains the initial ideal of the join:
$${\rm in}_\prec(I \ast J) \subseteq {\rm in}_\prec(I) \ast {\rm in}_\prec(J).$$
\end{lemma}

\begin{prop}\label{prop:compute}
Let $I$ be a radical ideal in a polynomial ring over an algebraically closed field.  Then
$$I^{(r)}  = I \ast m^r.$$
\end{prop}

\begin{proof}
We must compute the join $I \, \ast \,  m^r$.  By Lemma \ref{lem:intersect} it suffices to compute
$\cap_{\bfp \in V(I)} \, \,  (m_\bfp  \ast m^r)$.  The Zariski-Nagata theorem implies that it suffices to show that $m_\bfp \ast m^r = m_\bfp^r$, since then we will have deduced the equation
$$I \ast m^r  =   \left( \bigcap_{\bfp \in V(I)} m_\bfp \right)  \ast m^r  =      
\bigcap_{\bfp \in V(I)} (m_\bfp \ast m^r)  =  \bigcap_{\bfp \in V(I)} m_\bfp^r  = I^{(r)}.$$
So let $m_\bfp$ be a maximal ideal and let $\prec$ be any term order.  We have ${\rm in}_\prec(m_\bfp) = m$.    We compute $m \ast m^r$ by observing that if $\bfx^\bfa$ is a monomial it has expansion
$$(\bfy + \bfz)^\bfa =  \sum_{ \bfb + \bfc = \bfa}  \prod_{i =1}^n  {a_i \choose b_i} y_i^{b_i}  z_i^{c_i}.$$
The  monomial $\bfx^\bfa$ is in the join if and only if $(\bfy + \bfz)^\bfa \in m(\bfy) + m(\bfz)^r$.  Any term on the right hand side belongs to $m(\bfy)$ if any $b_i >0$.  If $\bfb = {\bf 0}$ then the term that appears is $\bfz^\bfa$, which belongs to $m(\bfz)^r$ if and only if $\bfx^\bfa \in m^r$.  Thus $m \ast m^r = m^r$.  Now ${\rm in}_\prec(m_\bfp \ast m^r)  \subseteq m \ast m^r =  m^r$.  On the other hand, any monomial $\bfx^\bfa \in m^r$ is the leading term of the monomial $(\bfx - \bfp)^\bfa \in m_\bfp^r$.  This implies that the polynomials of the form $(\bfx - \bfp)^\bfa$ with $|\bfa| = r$ form a Gr\"obner basis for $m_\bfp \ast m^r$ and hence $m_\bfp \ast m^r = m_\bfp^r$. 
\end{proof}

Proposition \ref{prop:compute} provides a straightforward proof of Proposition \ref{prop:secdiff}.

\smallskip

\noindent \emph{Proof of Proposition \ref{prop:secdiff}.}  Lemma \ref{lem:degreebound} implies that $I^{\{s-1\}} \subseteq m^s$.  The join operation respects inclusions, so
$$ I^{\{r+s-1\}} \, \, =  \, \,  I^{\{r\}} \ast I^{\{s-1\}} \, \,  \subseteq  \, \,  I^{\{r\}} \ast m^s  \, \,  = \, \,  \left(I^{\{r\}}\right)^{(s)}. $$
\qed


\section{Symbolic Powers of Edge Ideals}\label{sec:graph}

In this section, we explore properties of the symbolic powers of edge ideals.  In particular, we show that an edge ideal is $1$-differentially perfect if and only if the underlying graph is perfect.  As the secant ideals of edge ideals of perfect graphs were determined in \cite{Sturmfels2006}, this allows us to give explicit formulas for the symbolic powers of the edge ideals of perfect graphs.   The study of the symbolic powers of edge ideals was initiated in \cite{Simis1994} and further elaborated on in \cite{Bahiano2004, Gitler2006, Herzog2007}.  Our emphasis on perfect graphs, and the extensions to secant ideals of edge ideals in the next section, is new.  

Throughout this section and the next we use the fact that the generators of the symbolic powers of a squarefree monomial ideal do not depend on the underlying field.  Indeed, write $I = \cap P_i$ where $P_i$ is a monomial prime ideal.  Then $I^{(r)}  =  \cap P_i^r$ for all $r$.  Since this expression does not depend on the field, we will use the characteristic zero characterization of symbolic powers via derivatives  when computing symbolic powers of squarefree monomial ideals.

We first need some preliminary definitions regarding graphs and their  edge ideals.   
Let $G$ be an undirected graph with vertex set $[n] := \{1,2, \ldots, n\}$ and edge set $E(G)$.  The \emph{edge ideal} associated to $G$ is the ideal
$$I(G) \quad = \quad \left<  x_i x_j \, \, | \, \, \{i,j\} \in E(G) \right>.$$
This is a special case of the facet ideal construction, studied for example in \cite{Faridi2002}, however, the hypergraph language from \cite{Ha2006} will prove more useful.  Let $H$ be a hypergraph on  $[n]$ with edge set $E(H) = \{V_1, \ldots, V_d\}$.  We assume that $H$ is simple and loopless which means that $E$ contains no singletons and no pair of edges $V_i, V_j$ with $V_i \subseteq V_j$. 
If $H$ is a hypergraph its \emph{edge ideal} is the squarefree monomial ideal
$$I(H) \quad = \quad  \left<  m_V \, \, | \, \, V \in E(H) \right>$$
where $m_V =  \prod_{i \in V}  x_i$.   Clearly, every squarefree monomial ideal is an edge ideal for some $H$.  Edge ideals have been much studied in combinatorial commutative algebra \cite{Faridi2002, Ha2006, Jacques2004, Simis1994}.  The emphasis is often on homological properties of such ideals.

An \emph{independent set} of a hypergraph $H$ is a subset  $V \subset [n]$ such that no edge of $H$ is contained in $V$.  The \emph{independence number} $\alpha(H)$ of $H$ is the cardinality of the largest independent set of $H$.  A proper $k$-coloring of a hypergraph is an assignment of the numbers $\{1, 2, \ldots, k\}$ to the vertices of $H$ such that no edge of $H$ has all its vertices assigned the same number.  The \emph{chromatic number} $\chi(H)$ of $H$ is the smallest $k$ such that there exists a proper $k$-coloring of $H$.  Note that a proper $k$-coloring of $H$ is a partition of the vertices of $H$ into $k$ independent sets.  The notions of independence, coloring, and chromatic number generalize the usual notions for graphs.

A \emph{clique} of a graph $G$ is a collection of vertices of $G$ which form a complete subgraph.  The \emph{clique number} $\omega(G)$ is the cardinality of the largest clique in $G$.  Note that $\omega(G)$ is always a lower bound on the chromatic number $\chi(G)$.
 The complement of a graph $\overline{G}$ is the graph on $[n]$ such that $\{i,j\} \in E(\overline{G})$ if and only if $\{i,j\} \notin E(G)$.  Note that a clique of $G$ is an independent set of $\overline{G}$ and thus
\begin{equation}\label{eq:comp}
\alpha(G)  = \omega(\overline{G})  \quad \quad \mbox{ and } \quad \quad  \alpha(\overline{G}) = \omega(G).
\end{equation}
Furthermore, a proper $k$-coloring of $\overline{G}$ is a partition of $[n]$ into $k$ cliques of $G$.  If $V\subset [n]$, the induced subhypergraph $H_V$ is the hypergraph with vertex set $V$ and edges equal to the set of edges of $H$ that are contained in $V$.  If $G$ is a graph $G_V$ is called the induced subgraph.

\begin{defn}
A graph is \emph{perfect} if and only if $\omega(G_V) = \chi(G_V)$ for all $V \subseteq [n]$.
\end{defn}

Standard examples of perfect graphs are bipartite graphs and the comparability and incomparability graphs of partially ordered sets (posets).  Among the nonperfect graphs are the \emph{odd holes} which are odd cycles of length greater than three, and the complements of the odd holes.  The celebrated strong perfect graph theorem of Chudnovsky, Roberts, Seymour, and Thomas \cite{Chudnovsky2002} says that the odd holes and their complements are the only minimal imperfect graphs.  Though we will refer to minimal imperfect graphs throughout this section, we will not need to use this strong result about their structure.  Our main result in this section, in its simplest form, is the following theorem.

\begin{thm}\label{thm:perfect}
The edge ideal $I(G)$ is $1$-differentially perfect if and only if $G$ is perfect.
\end{thm}

In the special case where $G$ is a chordal graph, this result was already shown in \cite{Villarreal2006}.  Theorem \ref{thm:perfect} is remarkably similar to a theorem about the secant ideals of edge ideals.

\begin{thm} \cite[Prop. 3.4]{Sturmfels2006}
A graph is perfect if and only if all nonzero secant ideals $I(G)^{\{r\}}$ are generated in degree $ r+1$.  In this case $I(G)^{\{r\}}$ is generated by all monomials $m_V$ such that $V$ is a clique of $G$ of cardinality $r+1$.
\end{thm}

To prove Theorem \ref{thm:perfect} we will need a number of results that are interesting in their own right, and that will further strengthen and generalize the theorem.  A key idea in the proof is a combinatorial characterization of when a monomial belongs to the differential power of a squarefree monomial ideal, which is contained in Lemma \ref{lem:indep}.

\begin{defn}
Let $\bfa \in \nn^n$ be a nonnegative integer vector.  The \emph{blowup hypergraph} $H_\bfa$ is defined as follows.  The vertices of the blowup hypergraph are pairs $(i, k)$ such that $i \in [n]$ and $k \in [a_i]$.  A set of vertices  $(i_1, k_1), (i_2, k_2), \ldots, (i_m, k_m)$ forms an edge of $H_\bfa$ if and only if the $i_j$ are all distinct and $\{i_1, \ldots, i_m\}$ is an edge of $H.$
\end{defn}

Note that if $a_i = 0$ there are no vertices in the blowup hypergraph $H_\bfa$ with first index $i$.    In the case where $\bfa$ is a $0/1$ vector, $H_{\bfa}  = H_{{\rm supp}(\bfa)}$ where ${\rm supp}(\bfa) = \{i \, \, | \, \, a_i \neq 0 \}$.    In the case where $G$ is a graph, the blowup graph $G_\bfa$ is obtained from $G_{{\rm supp}(\bfa)}$ by replacing each vertex $i$ with a copy of the empty graph with $a_i$ vertices.

\begin{ex}
If $G$ is the graph consisting of a single edge and $\bfa = (a_1, a_2)$ with $a_1, a_2 > 0$, the blowup graph is the complete bipartite graph $G_\bfa = K_{a_1a_2}$.
\end{ex}

\begin{lemma}\label{lem:indep}
A monomial $\bfx^\bfa$   belongs to the symbolic power $I(H)^{(r)}$ if and only if 
\begin{equation}\label{eq:indep}
\deg(\bfx^\bfa) \geq r + \alpha(H_\bfa).
\end{equation} 
The monomial $\bfx^\bfa$ is a minimal generator of $I(H)^{(r)}$ if and only if the inequality (\ref{eq:indep}) is an equality and $\alpha(H_\bfa) = \alpha(H_{\bfa - {\bf e}_i})$ for all $i \in {\rm supp}(\bfa)$.
\end{lemma}

\begin{proof}
Let $\bfx^\bfa \in I(H)^{(r)}$.  This happens if and only if for every monomial $\bfx^\bfb$ with $\deg(\bfx^\bfb) \leq r-1$ dividing $\bfx^\bfa$, $\bfx^\bfa/ \bfx^\bfb \in I(H)$.  However, the monomial $\bfx^{\bfa-\bfb} \in I(H)$ if and only if $\bfx^{\bfa - \bfb}$ is divisible by $\bfx_\sigma$ for some edge $\sigma \in H$ if and only if the set of vertices 
$$\{(i, k) \in H_{\bfa} \, \, | \, \, i \in {\rm supp}(\bfa - \bfb), k \in [a_i - b_i] \}$$
is not an independent subset of $H_\bfa$.  This is guaranteed to happen for all $\bfx^\bfb$ with $\deg(\bfx^\bfb) \leq r-1$ if and only if $\alpha(H_\bfa) \leq \deg(\bfx^\bfa) - r$.

A monomial $\bfx^\bfa \in I(H)^{(r)}$ is a minimal generator if and only if $\bfx^{\bfa - {\bf e}_i} \notin I(H)^{(r)}$ for all $i \in {\rm supp}(\bfa)$.  This means that we want inequality (\ref{eq:indep}) to hold while
$$
\deg(\bfx^{\bfa - \bfe_i} ) < r + \alpha(H_{\bfa- \bfe_i}).
$$
As $\deg(\bfx^{\bfa - \bfe_i})  =  \deg( \bfx^\bfa) -1$  and $\alpha(H_{\bfa- \bfe_i}) \geq \alpha(H_\bfa) -1$, this can happen if and only if we have equality in (\ref{eq:indep}) and 
$\alpha(H_{\bfa- \bfe_i}) = \alpha(H_\bfa)$ for all $i$.
\end{proof}

A cover of a hypergraph $H$ is a subset $V \subset [n]$ such that every edge of $H$ contains at least one element of $V$.  The covering number $\tau(H)$ is the smallest cardinality of a cover of $H$.  Since the complement of any independent set is a cover, we have
$$\tau(H)  = n - \alpha(H).$$  Thus, Lemma \ref{lem:indep} can be reinterpreted in terms of covering numbers.

\begin{cor}
The generators of the $r$-th symbolic powers of $I(H)$ correspond to the blowup hypergraphs with covering number $\geq r$:
$$I(H)^{(r)} \, \, = \, \,  \left< \,  \bfx^\bfa \, \, | \, \,  \tau(H_\bfa)  \geq r  \, \right>.$$
A monomial $\bfx^\bfa$ is a minimal generator of $I(H)^{(r)}$ if and only if $H_\bfa$ has covering number $r$ but every induced hypergraph of $H_\bfa$ has covering number less than $r$.
\end{cor}

Again, the similarity with results about secant ideals of edge ideals is striking.

\begin{prop} \cite[Prop. 3.11]{Sturmfels2006} \label{prop:oldsecants}
The generators of the $r$-th secant ideal of $I(H)$ correspond to the induced subhypergraphs with chromatic number greater than $r$:
$$I(H)^{\{r\}}  =  \left<  m_V  \, \, | \, \, \chi(H_V) > r \right>.$$
A monomial $m_V$ is a minimal generator of $I(H)^{\{r\}}$ if and only if $H_V$ has chromatic number $r+1$, but every induced subhypergraph of $H_V$ has chromatic number less than $r+1$.
\end{prop}

Some useful facts about perfect graphs are summarized in the following proposition.  The last two of these are well-known results of Lovasz \cite{Lovasz1972}.

\begin{thm}\label{thm:lovacz}  Let $G$ be a perfect graph.
\begin{enumerate}
\item  Any induced subgraph of $G$ is a perfect graph.
\item The complement $\overline{G}$ is a perfect graph.
\item The graph $G_{v,k}$ obtained from $G$ by replacing the vertex $v$ with a complete graph of cardinality $k$ is a perfect graph.
\end{enumerate}
\end{thm}

Point 2 in Theorem \ref{thm:lovacz} is often called the Perfect Graph Theorem and Point 3 is known as the Replication Lemma.   Denote by $\mathcal{C}_k(G)$ the set of cliques of $G$ of cardinality greater than or equal to $k$.  We now have all tools in hand to prove the computational form of our main theorem in this section.

\begin{thm}\label{thm:perfectprod}
A graph $G$ is perfect  if and only if for all r 
$$I(G)^{(r)} =  \left<  \prod_{i = 1}^l  m_{V_i} \, \, | \, \,  V_i \in \mathcal{C}_2(G)  \mbox{ with }  \sum_{i =1}^l ( |V_i| -1) = r \right>. $$
\end{thm}

\begin{proof}
Call the ideal on the right-hand side of the equation $J_r(G)$.  First of all, notice that for any graph $J_r(G) \subseteq I(G)^{(r)}$.  One way to see this is to note that each clique $V$ gives a generator of $I(G)^{\{|V|-1\}}$ and thus, by the containment from (\ref{eq:diffperfect}), we deduce the desired containment.  For an alternate proof that does not reference secant ideals, simply take all partial derivatives of order $r-1$ of a monomial of the form $\prod m_{V_i}$  such that $\sum (|V_i| -1) = r$.  Such partial derivatives will either  be zero, or divisible by at least one edge $x_i x_j$ of $V$.

Suppose that $G$ is a perfect graph and let $\bfx^\bfa$ be a monomial in $I(G)^{(r)}$.  We wish to show that $\bfx^\bfa \in J_r(G)$.  We may assume, without loss of generality, that $\bfx^\bfa$ is a minimal generator of $I(G)^{(r)}$ and thus it satisfies $\deg(\bfx^\bfa) = r + \alpha(G_\bfa)$.  Since $G$ is perfect, so is the graph $\overline{G_\bfa}$ since it is obtained from $G$ by passing to an induced subgraph, taking the complement,  and replacing vertices of the resulting graph with complete graphs (this last part is the complementary operation to replacing vertices of $G_\bfa$ with empty graphs).  Since $\overline{G_\bfa}$ is perfect, it has a proper coloring using precisely $\alpha(G_\bfa) = \omega(\overline{G_\bfa})$ colors.  This coloring is a partition of the vertices of $\overline{G_\bfa}$ into $\alpha(G_\bfa)$ parts, each of which is a clique of $G_\bfa$.  Denote these cliques by $V_1, V_2, \ldots, V_{\alpha(G_\bfa)}$.  Now any clique $V = \{(i_1,k_1), \ldots (i_l, k_l) \} $ of the graph $G_\bfa$ maps to a clique $V' = \{i_1, \ldots, i_l\}$ of cardinality $l$ by deleting the second coordinate.  This coloring  of $G_\bfa$ yields the factorization
$$\bfx^\bfa  =  \prod_{i = 1}^l  m_{V'_i}.$$
We claim that the factorization on the right hand side of this equation implies that $\bfx^\bfa \in J_r(G)$.  To see why, we compute the sum:
$$
\sum_{i =1}^{\alpha(G_\bfa)}  (|V'_i| - 1)  \, \,   = \, \,     \sum_{i=1}^{\alpha(G_\bfa)} |V'_i|  - \alpha(G_\bfa) \, \,   =  \, \,    \deg(\bfx^\bfa)  - \alpha(G_\bfa) \, \,   = \, \,  r.
$$
We can remove all the cliques $V'_i$ of cardinality one without changing this sum.  The resulting monomial belongs to $J_r(G)$ and divides $\bfx^\bfa$.

Now suppose that $G$ is not a perfect graph.  We will show that there exists an $r$ such that $I(G)^{(r)} \neq J_r(G)$.  It suffices to consider the case where $G$ is a minimal imperfect graph (every subgraph of $G$ is perfect).  This implies that $\chi(\overline{G}) = \alpha(G) + 1$.  Let $x_{[n]}$ be the product of all indeterminates.  Note that $x_{[n]} \in I(G)^{(n - \alpha(G))}$ since 
$$\deg(x_{[n]}) \, \, = \, \, n \, \, = \, \,  (n- \alpha(G)) + \alpha(G).$$
We claim that $x_{[n]} \notin  J_{n - \alpha(G)}(G)$.  If it were, following the argument in the preceding paragraph in reverse, there would be a proper coloring of $\overline{G}$ using $\alpha(G) < \chi(\overline{G})$ colors.  This is a contradiction.
\end{proof}

\noindent  {\em  Proof of Theorem \ref{thm:perfect}.}
A monomial ideal $I$ is $1$-differentially perfect if and only if every monomial $\bfx^\bfa$ in the symbolic power $I^{(r)}$ can be written in the form:
$$\bfx^\bfa  =  \prod_{i = 1}^l  \bfx^{\bfb_i}$$
where each monomial $\bfx^{\bfb_i}  \in  I^{\{s_i\}}$ and such that the $s_i$ satisfy $\sum s_i  = r$. Now if $G$ is a perfect graph, $I(G)^{\{r\}}$ is generated by the cliques of cardinality $r+1$ in $G$.  Thus, by Theorem \ref{thm:perfectprod}, $I(G)$ is $1$-differentially perfect, since we have proven that the generating sets of $I(G)^{(r)}$ have the desired form.  

On the other hand, suppose that $G$ is not perfect.  Without loss of generality, we can take $G$ to be a minimal imperfect graph.  A theorem of Lovasz \cite{Lovasz1972b} says that $G$ has precisely $\alpha(G)\omega(G) + 1 = n$ vertices.  The monomial $x_{[n]}$ that is the product of all the variables belongs to $I(G)^{(n - \alpha(G))}$.  As every subgraph of $G$ is perfect, if we had
$$x_{[n]} \in  I(G)^{\{n - \alpha(G)\}}  +  \sum_{i = 1}^{n - \alpha(G) -1}  I(G)^{(i)} I^{(n - \alpha(G) - i)}$$
then either $x_{[n]} \in J_r(G)$ as defined in the proof of Theorem \ref{thm:perfectprod}, or $x_{[n]} \in I(G)^{\{n - \alpha(G)\}}$.  The first condition is impossible, as shown in the proof of Theorem \ref{thm:perfectprod} and the second condition could occur if and only if the chromatic number of $G$ was strictly greater than $ n - \alpha(G)$ by Proposition \ref{prop:oldsecants}.  However, the chromatic number of  a minimally imperfect graph is $\omega(G) +1$.  This leads to the inequality
$$\omega(G) + 1  >  \alpha(G)\omega(G)+ 1 - \alpha(G)$$
and thus
$$\frac{\omega(G)}{\omega(G) - 1}  >  \alpha(G).$$
As both $\omega(G)\geq 2$ and $\alpha(G) \geq 2$ for an imperfect graph, this is a contradiction. \qed
\smallskip

Theorem \ref{thm:perfectprod} implies a number of results about symbolic powers of edge ideals that appear in the literature.

\begin{cor}\cite{Simis1994}
An edge ideal $I(G)$ is normally torsion free if and only if $G$ is a bipartite graph. 
\end{cor}

\begin{proof}
Let $I(G)$ be a bipartite graph.  Since bipartite graphs are perfect, we know that $I(G)$ is $1$-differentially perfect.  As $I(G)^{\{2\}} = 0$ for any bipartite graph, we know that $I(G)^{(r)} = I(G)^{r}$ for all $r$.  On the other hand, if $G$ is not bipartite it must have an odd cycle $C$ of length $2r -1$.  The monomial $x_C \in I(G)^{(r)}$ but is not in $I(G)^r$.
\end{proof}

Lemma \ref{lem:indep} is an useful tool even when $G$ is not a perfect graph.  In particular, it allows us to explicitly characterize the minimal generators of $I(G)^{(r)}$ for small $r$.

\begin{cor}\label{thm:2power}
For any graph $G$,
$$I(G)^{(2)}  =  I(G)^{\{2\}} + I(G)^2.$$
In particular, $I(G)^{(2)}$ is generated by cubics of the form $x_{i}x_{j}x_k$ such that $\{i,j,k\}$ is a triangle in $G$ and quartics of the form $x_ix_jx_kx_l$ such that $\{i,j\}$ and $\{k,l\}$ are edges of $G$.  
\end{cor}

\begin{proof}
We already know the containment $I(G)^{\{2\}} + I(G)^2 \subset I(G)^{(2)}$.  Now suppose that $\bfx^\bfa$ is a minimal generator of $I(G)^{(2)}$ and let $G_\bfa$ be the blowup graph.  Lemma \ref{lem:indep} implies that the largest independent set of $G_\bfa$ has cardinality two less than the number of vertices.  Let $A$ denote such an independent set.  If $A$ has cardinality $1$, then $G_\bfa$ must be a triangle and hence $\bfx^\bfa \in I(G)^{\{2\}}$.  So suppose that the cardinality of $A > 1$.  Let $v_1$ and $v_2$ be the two vertices of $G_\bfa$ not in $A$.  These two vertices must each have an edge incident to $A$ and there must exist two disjoint vertices $w_1, w_2 \in A$ such that $v_1w_1$ and $v_2w_2$ are edges of $G_\bfa$.   Suppose that  $w_1$ and $w_2$ did not exist, that is both $v_1$ and $v_2$ were only incident to $w \in A$.   Then either $v_1v_2$ is an edge, in which case $wv_1v_2$ project to a triangle dividing $\bfx^\bfa$, or there is no edge between $v_1$ and $v_2$ in which case  $A \setminus \{w\} \cup \{v_1, v_2\}$ would be a larger independent set in $G_\bfa$.  But then $v_1w_1$ and $v_2w_2$ project to a pair of edges dividing $\bfx^\bfa$ and thus $\bfx^\bfa \in I(G)^{2}$.
\end{proof}

Despite the connection between secant ideals and symbolic powers that has driven many of the results in this section, Corollary \ref{thm:2power} shows that  the symbolic power will generally record much coarser information than the secant ideal.  Indeed, the symbolic square of an edge ideal is always generated in degrees three and four, whereas the secant of an edge ideal can require generators of arbitrarily large odd degree \cite[\S 3]{Sturmfels2006}.  The minimal generators of $I(G)^{\{2\}}$ of degree $\geq 5$ are all divisible by one of the quartics in $I(G)^2$ and so are ``lost'' when taking the symbolic square.


The characterization given for $1$-differentially perfect edge ideals can be extended to arbitrary ideals generated by quadratic monomials.  To do this, we need to replace the symbolic power with the differential power in a polynomial ring over a field of characteristic zero.  In the setting where $I$ is not a radical ideal,  $I$ is $r$-differentially perfect if
$$\left(I^{\{r\}}\right)^{<s>}  = I^{\{r+s-1\}} +  \sum_{i = 1}^{s-1}  \left(I^{\{r\}}\right)^{<i>} \left( I^{\{r\}} \right)^{<s-i>}$$
holds for all $s$.

\begin{thm}\label{thm:squares}
Let $I = I(G) +  \langle x_i^2 \, \,  | i \in \sigma \rangle$ be an ideal generated by quadratic monomials.  Then $I$ is $1$-differentially perfect if and only if $G$ is a perfect graph.  
\end{thm}

\begin{proof}
We define a new blowup graph $G^\sigma_\bfa$ which takes into account the square elements $x_i^2$ in $I$.  Namely, if $\bfa \in \nn^n$, $G_\bfa$ has vertices $(i,k)$ such that $i \in [n]$ and $k \in [a_i]$.  A pair of vertices $(i_1, k_1)$ $(i_2, k_2)$ forms an edge if $i_1i_2$ is an edge of $G$ or if $i_1 = i_2 \in \sigma$.  We claim that for any graph $G$, $\bfx^\bfa \in I^{<r>}$ if and only if $\deg(\bfx^\bfa) \geq r + \alpha(G^\sigma_\bfa)$.  The proof is the same as that of Lemma \ref{lem:indep}, except that if two vertices $(i, k_1)$ and $(i, k_2)$ are connected by an edge and remain after removing the vertices indexed by $\bfb$, then the monomial $\bfx^{\bfa - \bfb}$ is divisible by $x_i^2$.  

Now suppose that $G$ is a perfect graph.  First we need the characterization of the generators of $I^{\{r\}}$.  Such a characterization is implicit in \cite[Theorem 3.12]{Sturmfels2006}.  In particular, let $\bfr = (r,r,\ldots, r) \in \nn^n$ and let ${\bf 1} = (1,1,\ldots, 1)  \in \nn^n$.  Then the generators of $I^{\{r\}}$ are the $r+1$ element cliques in the blowup graph $G^\sigma_{\bfr + {\bf 1}}$, by the correspondence that cliques $(i_0,k_0), \ldots, (i_{r}, k_{r})$ correspond to monomials $x_{i_0} \cdots x_{i_{r}}$.  This correspondence allows us to simply follow the proofs of Theorems \ref{thm:perfectprod} and \ref{thm:perfect} to deduce that $I$ is $1$-differentially perfect.  

Conversely, if $G$ is a minimal imperfect graph, the argument in the proof of Theorem \ref{thm:perfect} shows that the monomial $x_{[n]}$ is a generator of $I^{<n -\alpha(G)>}$ but not in the ideal $$I^{\{n-\alpha(G)\}}  + \sum_{i =1}^{n - \alpha(G) -1}  I^{<i>} I^{<n - \alpha(G) - i >}.$$
Hence, $I$ is not $1$-differentially perfect.
\end{proof}


\section{Symbolic Powers of Antichain Ideals}

Among the perfect graphs are the incomparability graphs of partially ordered sets (posets).  This class of graphs proves to be an important special case for combinatorial commutative algebra as many initial ideals of combinatorially defined ideals are edge ideals of such incomparability graphs.
If $P$ is a partially ordered set with ground set $[n]$, associate the edge ideal 
$$J(P) = \left<  x_ix_j \, \, | \, \, \mbox{neither }   i \preceq j \mbox{ nor  }   j \preceq i \mbox{ in } P \right>. $$
Alternately, the ideal $J(P)$ is generated by the two element antichains of $P$.  As the incomparability graphs of posets are perfect (this is a classic corollary of Dilworth's Theorem), the generators of the secant ideals $J(P)^{\{r\}}$ are precisely the $r+1$ element antichains of $P$.  We call such secant ideals the antichain ideals of the poset $P$.  Denote the set of all antichains of $P$ of cardinality greater than or equal to $k$ by $\ca_k(P)$.
Thus we deduce:

\begin{cor}\label{cor:poset}
The symbolic powers of the poset ideal $J(P)$ are:
$$J(P)^{(s)}  =  \left< \prod_{i =1}^l  m_{A_i}  \, \, | \, \,   A_i \in \ca_2(P) \mbox{ with }  \sum_{i = 1}^l( |A_i| -1)  = s  \right>.$$
\end{cor}

Corollary \ref{cor:poset} has a far-reaching generalization to the symbolic powers of the antichain ideals $J(P)^{\{r\}}$.  The main result of this section will be the following theorem, characterizing the generating sets of the symbolic powers of the antichain ideals.

\begin{thm}\label{thm:poset}
The antichain ideals $J(P)$ are differentially perfect.  In particular, the symbolic powers of the antichain ideal $J(P)^{\{r\}}$ are: 
\begin{equation}\label{eq:thmposet}
\left( J(P)^{\{r\}}\right)^{(s)}  = \left<  \prod_{i =1}^l  m_{A_i}  \, \, | \, \,  A_i \in \ca_{r+1}(P) \mbox{ with }  \sum_{i =1}^l ( |A_i| - r )  = s \right>.
\end{equation}
\end{thm}

\begin{ex}\label{ex:badperfect}
It should be noted that the natural generalization of Theorem \ref{thm:poset} to arbitrary perfect graphs is false.  Indeed, consider the graph on six vertices that is the graph of the triangulation of a triangle, with edge set $E = \{12,13,23,24,25,35,36,45,56\}$.  This graph is easily seen to be perfect.  The secant square of the graph ideal $I(G)$ is generated by four cubics corresponding to the four triangles in $G$:
$$I(G)^{\{2\}}  = \left<  x_1x_2x_3, x_2x_3x_5, x_2x_4x_5, x_3x_5x_6 \right>.$$
The product of all the variables $x_1x_2x_3x_4x_5x_6$ is in the symbolic power $(I(G)^{\{2\}})^{(2)}$ but is not divisible by a clique of size four in $G$ (there are none) or the product of two cliques of size three. \qed
\end{ex}

The proof of Theorem \ref{thm:poset} depends in a crucial way on Greene's Duality Theorem for posets \cite{Greene1976}.  See \cite{Britz2001} for a recent survey of the duality theorem with many extensions, corollaries, and applications.   The duality theorem asserts a remarkable coincidence between two sequences of numbers associated to a poset.  For $i = 0, 1, 2, \ldots$ let $a_i$ (respectively $c_i$) be the maximal cardinality of the union of $i$ antichains (resp. chains) of $P$.  Define the sequences $\lambda_i$, $\overline{\lambda}_i$, by $\lambda_i =  a_i - a_{i-1}$ and $\overline{\lambda}_i = c_i - c_{i-1}$ for $i \geq 1$.

\begin{thm}[Duality Theorem for Finite Posets]
For any finite poset $P$, the sequences $\lambda = (\lambda_1, \lambda_2, \ldots)$ and $\overline{\lambda} = (\overline{\lambda}_1, \overline{\lambda}_2, \ldots)$ are non-increasing and form conjugate partitions of $n = |P|$.   
\end{thm}

Note that the graph in Example \ref{ex:badperfect} fails to satisfy the duality theorem (where antichain is replaced with clique and chain is replaced with independent set).  Thus, the obstruction to generalizing Theorem \ref{thm:poset} seems to be whether or not the duality theorem fails for a perfect graph $G$.  Indeed, as our proof will show, Theorem \ref{thm:poset} generalizes to any perfect graph $G$ with the property that all blowup graphs $G_\bfa$ satisfy the duality theorem with respect to cliques and independent sets.  This statement is summarized in Theorem \ref{thm:Greene}.

To prove Theorem \ref{thm:poset} we need to establish some basic facts about blowup hypergraphs in the context of incomparability graphs of posets as well as the relations to the partitions described by the duality theorem.  Given a poset $P$, the independent sets of the incomparability graph $J(P)$ are the  chains of $P$.  The hypergraph $H_r(P)$ such that the  antichain ideal $J(P)^{\{r\}}$ is the edge ideal of $H_r(P)$ has all $r+1$ element antichains of $P$ as its edges.  Thus, the independent sets of  $H_r(P)$ are all unions of $r$ chains of $P$.
This implies that the partition $\overline{\lambda}$ contains information about the sizes of independent sets in $H_r(P)$.

\begin{lemma}\label{lem:posetindep}
The cardinality of the largest independent set of $H_r(P)$ is the sum
$$\alpha(H_r(P)) =  c_r  = \sum_{i =1}^r  \overline{\lambda}_i.$$
\end{lemma}

Now let $\bfa \in \nn^n$.  To decide whether or not $\bfx^\bfa \in (J(P)^{\{r\}})^{(s)}$, we need to come to terms with the blowup hypergraph $H_r(P)_\bfa$.  These will turn out to be hypergraphs whose edges are antichains in related posets.

\begin{defn}
Let $P$ be a poset and $\bfa \in \nn^n$ a nonnegative integer vector.  The \emph{blowup poset} $P_\bfa$ is the new partially ordered set with ground set  consisting of all pairs $(i,j)$ such that $i \in P$ and $j \in [a_i]$ and subject to the ordering $(i,j) <  (k,l)$ if $i < k$ in $P$ or if $i = k$ and $j < l$.
\end{defn}

Note that if $a_i = 0$ there are no elements of $P_\bfa$ with first coordinate $i$.  The blowup poset $P_\bfa$ obtained from $P$ by replacing each element $i$ with a chain of length $a_i$.  

\begin{lemma}\label{lem:blowupposet}
$$H_r(P)_\bfa  =  H_r(P_\bfa).$$
\end{lemma}

\begin{proof}
A collection of elements $(i_0, j_0), \ldots, (i_{r}, j_{r})$ is an antichain of $P_\bfa$ if and only if $i_0, \ldots, i_{r}$ are distinct and form an antichain in $P$.  Since $H_r(P)_\bfa$ and $H_r(P_\bfa)$ have the same ground set, this implies that the edges of  $H_r(P)_\bfa$ and  $H_r(P_\bfa)$ are the same.
\end{proof}

\begin{lemma}
Let $\bfx^\bfa$ be a monomial and let $\lambda$ and $\overline{\lambda}$ be the partitions associated to the blowup poset $P_\bfa$.  Then
$$\deg(\bfx^\bfa)  =   s  +  \alpha( H_r(P_\bfa))$$
where $s  =  \sum_{i >r }  \overline{ \lambda}_i$.
\end{lemma}

\begin{proof}
We have 
$$
\deg(\bfx^\bfa)  \, \, =  \, \,  |P_\bfa|     \, \, =  \, \,   \sum_{i \geq 1}  \overline{\lambda}_i 
  \, \,  =  \, \,   \sum_{i > r} \overline{\lambda}_i  +  \sum_{i =1}^r  \overline{\lambda}_i 
\, \,  = \,\,   \sum_{i > r} \overline{\lambda}_i +  \alpha(H_r(P_\bfa))
$$
where the last equality follows from Lemma \ref{lem:posetindep}.
\end{proof}

\smallskip

\noindent  {\em  Proof of Theorem \ref{thm:poset}.}  
Let $K_{r,s}$ denote the monomial ideal on the right hand side of Equation \ref{eq:thmposet}.  First of all, note that $K_{r,s}  \subseteq  (J(P)^{\{r\}})^{(s)}$, since taking $s-1$ derivatives of any generating monomial of $K_{r,s}$ either gives zero or leaves at least one antichain of cardinality $\geq r+1$.  So our goal is to show the reverse containment  $(J(P)^{\{r\}})^{(s)} \subseteq K_{r,s}$.

Let $\bfx^\bfa \in (J(P)^{\{r\}})^{(s)}$.   We may suppose that the independence inequality for the degree is sharp, that is
$$\deg(\bfx^\bfa)  =  s +  \alpha(H_r(P_\bfa)).$$
Suppose that $k$ is the unique integer such that ${\lambda}_k \geq r+1$ while ${\lambda}_{k+1} \leq r$.  Let $P^*$ be the subposet of $P_\bfa$ whose elements consist of the union of any $k$ antichains yielding the maximal cardinality of the union, which is $\sum_{i =1}^k  {\lambda}_i$.  The new poset $P^*$ is the blowup poset $P_\bfb$ for a vector $\bfb$ such that $\bfx^\bfb$ divides $\bfx^\bfa$.  We will show that $\bfx^\bfb \in K_{r,s}$.

Associated to the new poset $P_\bfb$ are two new partitions $\lambda^*$ and $\overline{\lambda}^*$.
Since, by construction, $P_\bfb$ is the union of $k$ antichains, we have ${\lambda}^*_i  = 0$ for $i > k$.  Also, ${ \lambda}^*_k  \geq  {\lambda}_k \geq r+1$ since we must have the inequalities $a^*_i \leq a_i$ for all $i$ but $a^*_k  =  a_k$.  This in turn implies that $\overline{\lambda}^*_i  = k$ for all $i \in [r+1]$ and hence that,   
$$\sum_{i> r}  \overline{\lambda}^*_i \, \, = \, \,  \sum_{i > r}  \overline{\lambda}_i \, \, = \, \, s.$$
In particular, $\bfx^\bfb \in (J(P)^{\{r\}})^{(s)}$.  Let $A_1, \ldots, A_k$ be a partition of $P_\bfb$ into $k$ antichains.  Since ${\lambda}^*_k \geq r+1$, each of these antichains must have cardinality greater than or equal to $r+1$.  For each $i$, let $A'_i$ denote the projection of the antichain $A_i$ to $P$. We have
$$\bfx^\bfb  =  \prod_{i = 1}^k x_{A'_i}.$$
Now we evaluate the sum
$$\sum_{i =1}^k  (|A'_i| - r)  \, \, = \, \,  \sum_{i=1}^k |A'_i| -  rk  \, \, = \, \, \deg(\bfx^\bfb)  - rk  \, \, =  $$ 
$$  = \, \, \deg(\bfx^\bfb)  - \sum_{i = 1}^r  \overline{\lambda}^*_i  \, \, =  \, \, \deg(\bfx^\bfb) - \alpha(H_r(P_\bfb)) \, \, = \, \, s.$$
The third equality follows from the fact that $\overline{\lambda}^*_i  = k$ for all $1 \leq i \leq r+1$ and the fourth equality follows from Lemma \ref{lem:posetindep}.  This equation implies that  $\bfx^\bfb \in K_{r,s}$ and hence $(J(P)^{\{r\}})^{(s)} \subseteq K_{r,s}$.  \qed

\smallskip

In general, we can extend the proof of Theorem \ref{thm:poset} to edge ideals of graphs that satisfy Greene's Duality Theorem, with respect to the cliques and antichains.  Thus, to any graph, we define the sequence $a_i$ (respectively, $c_i$) for $i = 0,1,2, \ldots$ to be the maximal cardinality of the union of $i$ cliques (respectively, independent sets) of $G$.  The sequences $\lambda_i$ and $\overline{\lambda}_i$, are defined by $\lambda_i = a_i - a_{i -1}$ and $\overline{\lambda}_i = c_i - c_{i -1}$, for $i = 1,2, \ldots$.

\begin{defn}
A graph $G$ is called a \emph{Greene graph} if, for every vector $\bfa$, the blowup graph $G_\bfa$ has sequences $\lambda_i$, and $\overline{\lambda}_i$ that are nonincreasing and are dual partitions.
\end{defn}

As the proof of Theorem \ref{thm:poset} only depended on the fact that the incomparability graph of a poset is a Greene graph, we deduce:

\begin{thm}\label{thm:Greene}
If $G$ is a Greene graph, we have:
$$\left(I(G)^{\{r\}}\right)^{(s)}  \, \, = \, \,   \left< \prod_{i =1}^l  m_{V_i}  \, \, | \, \,  V_i \in \mathcal{C}_{r+1}(G) \mbox{ with }  \sum_{i =1}^l ( |V_i| - r )  = s \right>.$$
\end{thm}

\begin{ex}\label{ex:graphcomp}
Note that the converse to Theorem \ref{thm:Greene} does not hold.  In particular, consider the graph $G$ on six vertices with edge set $E(G) = \{14, 25,36,45,46,56\}$.  This graph is not a Greene graph because the sequence $\lambda = (3,1,2,0,...)$ is not a partition.  On the other hand, $G$ is perfect so $I(G)^{(s)}$ is generated  by the product of cliques for all $s$.  Furthermore $I(G)^{\{2\}} = \left< x_4x_5x_6 \right>$ and $I^{\{r\}} = \left<0\right>$ for all $r > 2$.  Thus, the symbolic powers of the secant ideals of $I(G)$ satisfy the conclusion of Theorem \ref{thm:Greene}.
\end{ex}

Lemma \ref{lem:blowupposet} together with Greene's Theorem imply that the incomparability graphs of posets are Greene graphs (hence, the name).  It is easy to see that the comparability graphs of posets are also Greene graphs, which will prove useful in Section 6.  Recall that such a comparability graph has as vertices the elements of the poset $P$, and $ij$ is an edge  if and only if either $i < j$ or $j < i$ in $P$.

\begin{prop}\label{prop:compgreene}
The comparability and incomparability graphs of a poset are Greene graphs.
\end{prop}

\begin{proof}
That the incomparability graph of a poset is a Greene graph is the content of Lemma \ref{lem:blowupposet}.  We must prove that the comparability graph of a poset is a Greene graph.
Let $P$ be the underlying poset, and $G$ the comparability graph of $P$.  It suffices to show that the blowup graph $G_\bfa$ is the comparability graph of an associated poset $P^\bfa$.  In this case, Greene's Theorem will imply that all the blowup graphs are Greene graphs.  Define a poset $P^\bfa$ as follows.  The elements of $P^\bfa$ are pairs $(i,j)$ such that $i \in [n]$ and $j \in [a_i]$.  We have a relation $(i,j) < (k,l)$ if and only if $i < j$ in $P$.  Thus the poset $P^\bfa$ is obtained from $P$ by replacing the element $i$ with an antichain of cardinality $a_i$.  A set of elements $(i_1, j_1), (i_1,j_2)$ forms an edge of the comparability graph $P^\bfa$ if and only if $(i_1, j_1), (i_1,j_2)$ are comparable in $P^\bfa$ if and only if $i_1$ and $i_2$ are comparable in $P$ if and only if $i_1i_2$ is an edge of $G_\bfa$.  Thus $G_\bfa$ is the comparability graph of $P^\bfa$.
\end{proof}

It is worth noting that every Greene graph is perfect, but not every perfect graph is a Greene graph.  In particular, the graphs from Examples \ref{ex:badperfect} and \ref{ex:graphcomp} and are not Greene graphs.  This class of graphs seem not to have been studied in the graph theory literature and so it is an interesting open problem to find a characterization of this subclass of perfect graphs.

\begin{ques}
\begin{enumerate}
\item  Is it sufficient to only check induced subgraphs in the definition of a Greene graph?  In other words, is there a replication lemma for the set of graphs that satisfy the duality theorem for all induced subgraphs?
\item  What collection of excluded induced subgraphs characterize Greene graphs?
\end{enumerate}
\end{ques}


\section{Delightful Term Orders}\label{sec:delight}

Besides the interesting combinatorial questions that  arise, one motivation for studying the symbolic powers of monomial ideals is to try to use this information to prove theorems about symbolic powers of general ideals.  This is because of the following proposition.

\begin{prop}\label{prop:include}
Let $\prec$ be a term order such that both $I$ and ${\rm in}_\prec(I)$ are radical ideals.  Then
$${\rm in}_\prec(I^{(r)})  \subseteq   {\rm in}_\prec(I)^{(r)}.$$
\end{prop}

\begin{proof}
Since both $I$ and ${\rm in}_\prec(I)$ are radical, we can use Proposition \ref{prop:compute} to compute the symbolic powers.  Indeed, we have
$${\rm in}_\prec(I^{(r)})  =  {\rm in}_\prec(I \ast m^r)  \subseteq {\rm in}_\prec(I)  \ast {\rm in}_\prec(m^r)  =  {\rm in}_\prec(I) \ast m^r  =  {\rm in}_\prec(I)^{(r)}.$$
The first and last equality follow from Proposition \ref{prop:compute}, the containment follows from Lemma \ref{lem:intersect}, and the middle equality follows because $m^r$ is a monomial ideal. 
\end{proof}

Thus, a strategy for constructing Gr\"obner bases (and hence generating sets) for the symbolic powers $I^{(r)}$ would be the following:
\begin{enumerate}
\item  Compute ${\rm in}_\prec (I)$ and give a combinatorial description for its minimal generators.
\item  Determine a combinatorial description of the symbolic power $ {\rm in}_\prec (I)^{(r)} $.
\item  Find a collection of polynomials $\mathcal{G} \subset I^{(r)}$ such that  $\left< {\rm in}_\prec(\mathcal{G}) \right> =  {\rm in}_\prec (I)^{(r)} $.
\item Deduce that $\mathcal{G}$ is a Gr\"obner basis for $I^{(r)}$.
\end{enumerate}

In this section, we explain how to pursue this strategy for some classic ideals of combinatorial commutative algebra, in particular, for determinantal and Pfaffian ideals.  Note that this is the same strategy that was described for computing secant ideals combinatorially in \cite{Sturmfels2006}.   In fact, there is a close connection between applying this method for secant ideals and for differential powers.  Recall the following definition for secants of ideals.

\begin{defn}
A term order $\prec$ is called $r$-\emph{delightful} for $I$ if 
$${\rm in}_\prec(I^{\{r\}}) =  {\rm in}_\prec(I)^{\{r\}}.$$
A term order $\prec$ is \emph{delightful} if it is $r$-delightful for all $r$.
\end{defn}

\begin{thm}\label{thm:gb}
Let $I$ and ${\rm in}_\prec(I)$ be radical and
suppose that $\prec$ is an $r$-delightful term order for $I$ for all $r \leq t$, and that for all $r \leq t$ and $s \leq u$, ${\rm in}_\prec(I)$ satisfies
\begin{equation}\label{eq:in}
\left({\rm in}_\prec(I)^{\{r\}}\right)^{(s)} \quad = \quad {\rm in}_\prec(I)^{\{r+s-1\}} + \sum_{i =1}^{s-1} \left({\rm in}_\prec(I)^{\{r\}}\right)^{(i)} \left({\rm in}_\prec(I)^{\{r\}}\right)^{(s-i)}. \end{equation}
For $r \leq t$ let $\cg_r = \{g_1^r, g_2^r, \ldots \}$ be a Gr\"obner basis for $I^{\{r\}}$ with respect to $\prec$.  Then for all $r \leq t$ and $s \leq \min(u, t - r + 1)$ the set of polynomials
$$\cg_{r,s} \quad  =  \quad  \left\{ \prod_{i =1}^l  g^{r_i}_{j_i}  \, \, | \, \, r_i \geq r, \, \, \sum_{i =1}^l  (r - r_i +1) = s \right\}$$
is a Gr\"obner basis for $(I^{\{r\}})^{(s)}$ with respect to $\prec$.  In particular, for all $r \leq t$ and $s \leq  \min(u, t-r+1)$, $I$ satisfies
$$\left(I^{\{r\}} \right)^{(s)}  \quad = \quad  I^{\{r+s -1\}} + \sum_{i=1}^{s-1}  \left(I^{\{r\}} \right)^{(i)} \left(I^{\{r\}} \right)^{(s-i)}.$$ 
\end{thm}

\begin{proof}
Since the initial ideal ${\rm in}_\prec(I)$ satisfies Equation \ref{eq:in},  the minimal generators of \linebreak $({\rm in}_\prec(I)^{\{r\}})^{(s)}$ have the form $\prod_{i=1}^l m^{r_i}_{j_i}$ such that $\sum_{i =1}^l(r_i - r +1 ) = s$ where $m^{r_i}_{j_i} \in {\rm in}_\prec(I)^{\{r_i \}}$.  However, since $r_i \leq t$  and $\prec$ is $r_i$ delightful, $m^{r_i}_{j_i}$ is the leading term of a polynomial in $f^{r_i}_{j_i} \in \cg_{r_i}$.  This implies that each monomial generator of $({\rm in}_\prec(I)^{\{r\}})^{(s)}$  is the leading term of a polynomial in $\cg_{r,s}$.   Since $\cg_{r,s} \subseteq (I^{\{r\}})^{(s)}$ by the containment (\ref{eq:diffperfect}) we deduce that the initial terms of $\cg_{r,s}$ generate the initial ideal ${\rm in}_\prec((I^{\{r\}})^{(s)} )$ by Proposition \ref{prop:include} and that $\cg_{r,s}$ is a Gr\"obner basis for $(I^{\{r\}})^{(s)}$.
\end{proof}

Sending $t$ and $u$  to infinity, we deduce the following corollary. 

\begin{cor}\label{cor:delightperfect}
Suppose that $I$ and ${\rm in}_\prec(I)$ are radical, that $\prec$ is a delightful term order for $I$, and that ${\rm in}_\prec(I)$ is differentially perfect.  Then $I$ is differentially perfect.
\end{cor}

Thus nice descriptions of the Gr\"obner bases of secant ideals and symbolic powers  seem to go hand-in-hand.   To conclude this section, we show how our combinatorial techniques can be used to derive Gr\"obner bases for the symbolic powers of some classical ideals.
Our first example concerns the ideals of minors of a generic matrix.

\begin{thm}\label{thm:generic}
The ideal $I_{mn}$, generated by the $2 \times 2$ minors of a generic $m \times n$ matrix $X_{mn}$, is differentially perfect.
\end{thm}

\begin{proof}
Let $\prec$ be any diagonal term order, that is, any term order that selects the main diagonal of any subdeterminant  of $X_{mn}$ as the leading term.  The $2 \times 2$ minors of $X_{mn}$ form a Gr\"obner basis of $I_{mn}$ with respect to $\prec$, and the initial ideal ${\rm in}_\prec(I_{mn}) = J(P_{mn})$ for the poset $[m] \times [n]$ subject to the ordering $(i,j) \leq (k,l)$ if and only if $i \geq k$ and $j \leq l$.  Thus the initial ideal ${\rm in}_\prec(I)$ is radical and differentially perfect.
Diagonal term orders are also delightful for $I_{mn}$ \cite[\S 4]{Sturmfels2006}.  Thus, by Corollary \ref{cor:delightperfect}, $I_{mn}$ is differentially perfect.  
\end{proof}

\begin{defn}
For any matrix $X$, $\mathcal{M}_r(X)$ is the union of the set of all $t \times t$ minors of $X$ for all $t \geq r$.  For a skew-symmetric matrix $Y$, $\mathcal{P}_r(Y)$ is the union of the set of all $2t \times 2t$ subPfaffians of $Y$ for all $t \geq r$.
\end{defn}

\begin{cor}
The $s$-th symbolic power  of the ideal $I^{\{r\}}_{mn}$ of $(r+1) \times (r+1)$ minors of a generic $m\times n$ matrix $X_{mn}$ are generated by products of minors.  In particular:
$$\left(I^{\{r\}}_{mn}\right)^{(s)} =  \left<  \prod_{i =1}^l  f_i  \, \, | \, \, f_i \in \mathcal{M}_{r+1}(X_{mn}),  \,\, \sum_{i =1}^l  (\deg f_i - r)  = s  \right>$$
and these products of minors form a Gr\"obner basis for $\left(I^{\{r\}}_{mn}\right)^{(s)}$ with respect to any diagonal term order.
\end{cor}

The usual diagonal term orders for symmetric minors and Pfaffians were shown to be delightful in \cite[\S 4]{Sturmfels2006}.  In both cases, the initial ideal for the $k= 2$ case (i.e. $2 \times 2$ minors and $4\times 4$ Pfaffians, respectively) is an antichain ideal $J(P)$.  Thus, by Theorem \ref{thm:gb}, the symbolic powers have Gr\"obner bases consisting of the obvious products of minors and Pfaffians, respectively.  We state these results in Theorems \ref{thm:symmetric} and \ref{thm:pfaffian}.

\begin{thm}\label{thm:symmetric}
The ideal $I_m$ of $2\times 2$ minors of a generic symmetric $m \times m$ matrix $X_m$ is differentially perfect.  In particular the set
$$\left \{\prod_{i =1}^l  f_i  \, \, | \, \, f_i \in \mathcal{M}_{r+1}(X_m), \, \, \sum_{i =1}^l  (\deg f_i -r) = s \right\}
$$
forms a Gr\"obner basis for the symbolic power $\left(I_m^{\{r\}}\right)^{(s)}$ with respect to any diagonal term order.
\end{thm}

In Theorem \ref{thm:symmetric}, the poset $P_{m}$ such that $J(P_m)  = {\rm in}_\prec(I_m)$ consists of all pairs $(i,j) \in [m] \times [m]$ such that $i \leq j$, subject to the ordering $(i,j) \leq (k,l)$ if and only if $i \geq k$ and $j \leq l$.

\begin{thm}\label{thm:pfaffian}
The ideal $I_m$ of $4 \times 4$ Pfaffians of a generic $m \times m$ skew-symmetric matrix $Y_m$ is differentially perfect.  In particular the set
$$\left \{\prod_{i =1}^l  f_i  \, \, | \, \, f_i \in \mathcal{P}_{r+1}(Y_m), \, \, \sum_{i =1}^l  (\deg f_i -r) = s \right\}
$$
forms a Gr\"obner basis for the symbolic power $\left(I_m^{\{r\}}\right)^{(s)}$ with respect to any antidiagonal term order.
\end{thm}

In Theorem \ref{thm:pfaffian}, the poset $P_m$ such that $J(P_m) = {\rm in}_\prec(I_m)$ consists of all pairs $(i,j) \in [m] \times [m]$ such that $i < j$, subject to the ordering $(i,j) \leq (k,l)$ if and only if $i \geq k$ and $j \geq l$.

\begin{rmk}
The arguments presented in the previous theorems work for ladder determinantal and Pfaffian ideals as well, which was originally treated in \cite{Bruns1998}.  Indeed, the diagonal (respectively, antidiagonal) term order is easily shown to be delightful in these cases and the initial ideal is still a poset ideal $J(P)$ for a modified poset.   Corollary \ref{cor:delightperfect} applies in the usual way.
\end{rmk}


\section{Symbolic Powers of Some Segre-Veronese Ideals}\label{sec:app}

Given vectors of nonnegative integers $\bfn = (n_1, \ldots, n_m)$ and $\bfd = (d_1, \ldots, d_m)$, the Segre-Veronese variety is the variety
$$V_{\bfn,\bfd}   \quad =  \quad  \nu_{d_1}( \pp^{n_1}) \times \cdots \times \nu_{d_m}(\pp^{n_m})  \subset  \pp^{N}$$
where $N  =  \prod  {d_i + n_i  \choose n_i} - 1$, $\nu_{d_i}$ denotes the $d_i$-uple Veronese embedding, and the products $\times$ denote the usual Segre product.  Segre-Veronese varieties naturally generalize the Segre varieties and Veronese varieties, and their secant varieties and symbolic powers pose many interesting questions \cite{Catalisano2005}.

The ideals of  Segre-Veronese varieties are, in many cases, the ideals of $2 \times 2$ minors of certain matrices, such as generic matrices, symmetric matrices, Hankel matrices, and catalecticant matrices.  In some cases, the ideals of the secant varieties are also generated by minors of matrices, though it seems difficult to characterize precisely when this happens.  In this section, we explore three cases where this occurs, showing the results by producing delightful term orders where the minors form Gr\"obner bases, and using the edge ideal structure of the initial ideals to realize Gr\"obner bases of the symbolic powers.

One way to view our approach to computing Gr\"obner bases of the symbolic powers of these combinatorial ideals and the ideals in Section \ref{sec:delight}, is that we are replacing the Knuth-Robinson-Schensted correspondence, used  in the standard proofs, with Greene's Duality Theorem.  Each of the proofs of the Gr\"obner basis results in these cases depends on finding a different straightening law which often uses the KRS correspondence.  It should be noted that our combinatorial approach to secants and symbolic powers is not entirely separate from the KRS correspondence.  Indeed, as shown in \cite{Britz2001}, the KRS correspondence is a corollary of the duality theorem.  Thus, our approach seems to extract the ``combinatorial essence'' of the problem and gives another explanation for why the KRS algorithm works.  Generally, we expect a KRS based approach to be successful for studying secant varieties and symbolic powers when there is an initial ideal that is the antichain ideal of  a \emph{wonderful poset} (see \cite{DeConcini1982}).   Among the examples in this section are ideals whose initial ideals are edge ideals of Greene graphs that are not incomparability graphs, where the KRS approach seems not to apply.


\subsection{The Rational Normal Curve $\nu_d(\pp^1)$}\label{sec:hankel}

Our first example concerns the secants and symbolic powers of the ideal $I_d$ of the rational normal curve $\nu_d(\pp^1)$ embedded in $\pp^d$ in the standard toric embedding.  The ideal $I_d$ is generated by the $2 \times 2$ minors of the $2 \times d$ Hankel matrix:
$$X_1  =  \begin{pmatrix}
x_0 & x_1 & x_2 &  \cdots & x_{d-1} \\
x_1 & x_2 & x_3 & \cdots & x_d
\end{pmatrix}.$$
The secant ideals $I^{\{r\}}_d$ are generated by the $(r+1) \times (r+1)$ minors of the  $(r+1) \times (d-r+1)$ Hankel matrix:
$$X_{r}  = \begin{pmatrix}
x_0 & x_1 & x_2 &  \cdots & x_{d-r} \\
x_1 & x_2 & x_3 & \cdots &   x_{d-r+1} \\
 \vdots & \vdots & \vdots & \ddots & \vdots \\
 x_r & x_{r+1} & x_{r+2}  &  \cdots & x_d  \end{pmatrix}.$$
A standard monomial theory for the minors of Hankel matrices was developed by Conca \cite{Conca1998} to show that: the $(r +1) \times (r+1)$ minors are Gr\"obner bases of the secant ideals and that the appropriate products of minors form Gr\"obner bases for the symbolic powers of the secant ideals.  A straightening law is also developed for minors of Hankel matrices to give a primary decomposition and description of the initial ideals of the ordinary powers of the secant ideals.  We will show how to derive the first two of these results from our combinatorial framework, together with results concerning the connections between Gr\"obner bases of toric ideals and triangulations of polytopes \cite{Sturmfels1996}.

Let $Z_d$ denote the \emph{zigzag poset}, whose elements are the numbers $\{0,1,\ldots, d \}$ and whose only relations are $2i -1 <  2i$ and $2i > 2i+1$ for all $i$.  The zigzag poset $Z_8$ is pictured in the figure.

\begin{figure}[h]
\begin{center}\includegraphics{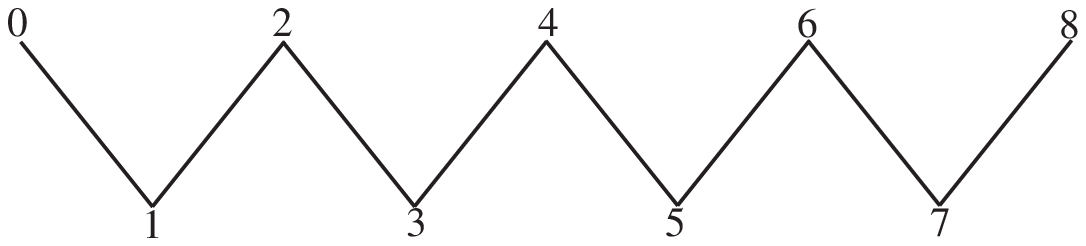} \end{center}
\end{figure}

Let $\prec$ be any term order that selects the main diagonal of every maximal minor of the Hankel matrices $X_{r}$ as leading term.  One such term order is the lexicographic order with $x_0 \succ x_1 \succ \cdots \succ x_d$.

\begin{prop}
The initial ideal of the ideal $I_d$ of $2 \times 2$ minors of the generic Hankel matrix $X_1$ with respect to any diagonal term order $\prec$ is the antichain ideal $J(Z_d)$.  The $2 \times 2$ minors of $X_1$ form a Gr\"obner basis for $I_d$.
\end{prop}

\begin{proof}
The rational normal curve is a toric variety and the associated polytope is a segment of length $d$.  The diagonal term order induces a triangulation of this segment into $d$ segments of unit length.  The resulting simplicial complex has facets $\{0,1\}, \{1,2\}, \ldots,$ $\{d-1, d\}$.  As this triangulation is unimodular and the minimal nonfaces of the associated simplicial complex are all edges, the initial ideal ${\rm in}_\prec(I_d)$ is generated by squarefree quadrics, which have the form $x_ix_j$ such that $j > i+1$.  Each such quadric is the diagonal leading term of the $2 \times 2$ minor of the submatrix:
$$\begin{pmatrix}
x_i & x_{j-1} \\
x_{i+1} & x_j
\end{pmatrix}$$
and is also an incomparable pair in the zigzag poset $Z_d$.  Conversely, every such incomparability pair is the leading term of such a $2 \times 2$ minor.
\end{proof}

\begin{thm}
Any diagonal term order $\prec$ is delightful for the ideal $I_d$ of the $2\times 2$ minors of the generic Hankel matrix $X_1$.  The $(r+1) \times (r+1)$ minors of $X_{r}$ generate $I_d^{\{r\}}$ and form a Gr\"obner basis with respect to $\prec$.
\end{thm}

\begin{proof}
The initial ideal ${\rm in}_\prec(I_d)$ is the antichain ideal $J(Z_d)$ of the zigzag poset.  The antichains of $Z_d$ consist of sequences $i_0, i_1, \ldots, i_{r}$ such that $i_{k+1} > i_k + 1$.  Each such antichain is the diagonal leading term of an $(r +1) \times (r + 1)$ minors of $X_{r+1}$ of the form:
$$\begin{pmatrix}
x_{i_0} & x_{i_1 -1} & x_{i_2 - 2} & \cdots & x_{i_r -r} \\
x_{i_0 +1} & x_{i_1} & x_{i_2 -1} & \cdots &  x_{i_r - r + 1} \\
\vdots &  \vdots &  \vdots & \ddots & \vdots \\
x_{i_0 + r} & x_{i_1 + r-1}  & x_{i_2 + r -2} &  \cdots & x_{i_r} 
\end{pmatrix}.
$$
Each such minor belongs to the secant ideal $I_d^{\{r\}}$ by elementary linear algebra.  This implies that ${\rm in}_\prec(I_d)^{\{r\}} =  {\rm in}_\prec(I_d^{\{r\}})$ which implies that $\prec$ is delightful.
\end{proof}

\begin{cor}\label{cor:hankel}
The ideal $I_d$  of $2 \times 2$ minors of the generic Hankel matrix $X_1$ is differentially perfect.
In particular, the set 
$$\left\{ \prod_{i =1}^l f_{j_i}^{r_i} \, \, | \, \, f^{r_i}_{j_i} \in \mathcal{M}_{r_i + 1} (X_{r_i}), \, \, r_i \geq r,   \, \,  \sum_{i =1}^l (\deg f^{r_i}_{j_i} - r) = s \right\}$$ 
 forms a Gr\"obner basis for the symbolic power $(I_d^{\{r\}})^{(s)}$ with respect to any diagonal term order.
\end{cor}

\begin{proof}
The diagonal term order $\prec$ is delightful for $I_d$.  The initial ideal ${\rm in}_\prec(I_d)$ are poset ideals $J(Z_d)$, and hence also edge ideals for a Greene graph $G_d$, the incomparability graph of $Z_d$.  The result follows by Theorem \ref{thm:gb} and Corollary \ref{cor:delightperfect}.
\end{proof}


\subsection{ The Surface $\nu_d(\pp^1) \times \nu_2(\pp^1)$}\label{sec:symhankel}

Let $\kk[x] := \kk[x_{ij} \, \, | \,\, i = 0,1, ..., d, j = 0,1,2]$ be the polynomial ring in $3(d+1)$ indeterminates and let $X_k$ denote the block Hankel matrix:

$$X_k  = \left( \begin{array}{ccccc}
A_0 & A_1 & A_2 & \cdots & A_{d-k} \\
A_1 & A_2 & A_3 & \cdots & A_{d-k+1} \\
\vdots & \vdots & \vdots & \ddots & \vdots \\
A_k & A_{k+1} & A_{k+2} & \cdots & A_d  \end{array} \right)$$
where each $A_i$ is a $2 \times 2$ matrix of indeterminates:
$$A_i =  \begin{pmatrix}
x_{i0} & x_{i1} \\  x_{i1} & x_{i2} \end{pmatrix}. $$

Let $I_d$ be the ideal generated by the $2 \times 2$ minors of $X_{\lfloor d/2 \rfloor}$.  Let $\prec$ be any term order that selects the main diagonal as the leading term of any minor of any of the matrices $X_k$.  One such term order is the lexicographic order with $x_{i_1j_1} \prec x_{i_2j_2}$ if $i_1 > i_2$ or if $i_1 = i_2$ and $j_1 > j_2$.  First of all, we claim that these $2 \times 2$ minors form a Gr\"obner basis for $I_d$, that they generate the ideal $I( \nu_d(\pp^1) \times \nu_2(\pp^1))$,  and that the initial ideal is an antichain ideal for a poset $P_d$.

\begin{figure}[h]
\begin{center}\includegraphics{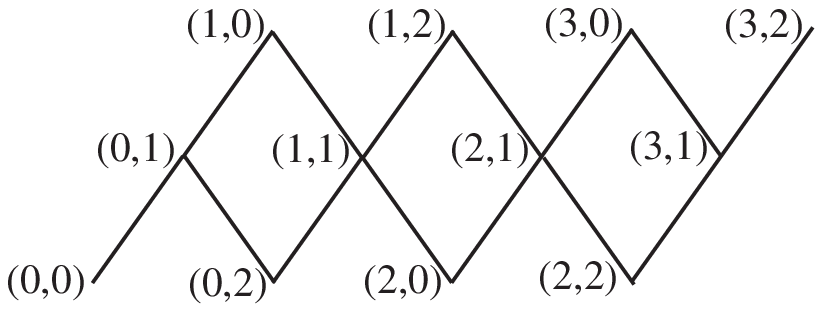} \label{fig:poset}
\end{center}
\end{figure}

In particular, let $P_d$ be the poset on the pairs $(i,j)$ subject to the following covering relations:
$$(2i,0) \prec (2i,1), \quad (2i,2) \prec (2i, 1), \quad (2i+1, 1) \prec (2i+1,0), \quad (2i+1, 1) \prec (2i+1,2), $$
$$(2i,1) \prec (2i-1,2), \quad  (2i,1) \prec (2i+1, 0), \quad (2i,2) \prec (2i+1, 1), \quad (2i,0) \prec(2i-1,1).$$
The poset $P_3$ is pictured in Figure \ref{fig:poset} and the basic pattern continues for larger $d$.

\begin{prop}
The $2 \times 2$ minors of $X_{\lfloor d/2 \rfloor}$ form a Gr\"obner basis for the ideal $I_d$ with respect to any diagonal term ordering.  The initial ideal of $I_d$ is the antichain ideal $J(P_d)$. 
\end{prop}

\begin{proof}
To show these results, recall that any Segre-Veronese variety is a toric variety, and thus the vanishing ideal $I_d$ is a toric ideal.  Let $\phi_d$ be the ring homomorphism:
$$\phi_d :  \kk[x]  \rightarrow  \kk[t,u,v]$$
$$x_{ij}  \mapsto  tu^iv^j.$$
Let $J_d = \ker \phi_d$ be the toric ideal that defines this Segre-Veronese variety $ \nu_d(\pp^1) \times \nu_2(\pp^1)$.  The relations in any toric ideal are determined by the combinatorics of the associated configuration of exponent vectors appearing in the parametrization.  In our case, this consists of the vectors $(1,i,j) \in \nn^3$ where $i \in \{0,1,\ldots, d\}$ and $j \in \{0,1,2\}$.  Since this vector configuration is homogenous (all the points lie on a plane that does not pass through the origin), we can reduce to a 2-dimensional configuration of points.  In our case, these are the $3(d+1)$ integer points in the rectangle $[0,2] \times [0,d]$.

To construct a quadratic initial ideal, we use the fact that the initial complexes of toric ideals are the regular triangulations of the corresponding point configurations (see \cite{Sturmfels1996} for background).  In particular, the triangulation with respect to the lexicographic term order described above, is depicted in Figure \ref{fig:triangulate}.  The pattern of the triangulation continues to the right with increasing $d$.  

\begin{figure}[h]
\begin{center}
\includegraphics{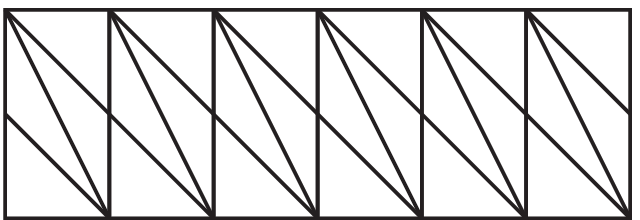} \label{fig:triangulate}
\end{center}
\end{figure}

As the minimal nonfaces of the triangulation are all edges and each triangle has area $\frac{1}{2}$, the corresponding initial ideal is squarefree and generated by quadrics that correspond to the non-edges in the triangulation.  First we will show that these nonedges are precisely the incomparable pairs in the poset $P_d$.  This is equivalent to showing that every edge in the triangulation is a comparable pair in the poset.  There are seven different types of edges in the triangulation, namely, the edges 
$$((i,0),(i+1,0)), ((i,0), (i,1)), ( (i,1),(i+1, 0)), ((i,1), (i,2))$$
$$((i,2),(i+1,0)), ((i,2),(i+1,2)), ((i,2),(i+1,2)).$$
Each of these corresponds to a comparable pair in the poset, since for example $((i,0),(i+1,0))$ is a comparable pair in $P_d$ for any $i$.  Conversely, every lexicographically ordered comparable pair falls into one of the seven classes.

Now we claim that each of these nonedges is the initial term of a $2 \times 2$ minor of the matrix $X_{\lfloor d/2 \rfloor}$ with respect to the diagonal term ordering.   This is an immediate consequence of Lemma \ref{lem:lowerdiag} below.
\end{proof}

Given an arbitrary $m \times n$ matrix $Y$, the $k$-th superdiagonal consists of all entries of the form $y_{i+k,i}$ of $X$.  Given a pair of indices $i \in [m] $ and $j \in [n]$ the $ij$ lower right submatrix $Y^{\searrow}_{ij}$ is the submatrix
$$Y^{\searrow}_{ij} = \begin{pmatrix}
y_{i,j} & y_{i,j+1} & \cdots & y_{i,n} \\
y_{i+1,j} & y_{i+1,j+1} & \cdots & y_{i+1,n} \\
\vdots & \vdots & \ddots &  \vdots \\
y_{m,j} & y_{m,j+1} & \cdots & y_{mn} \end{pmatrix}.
$$

\begin{lemma}\label{lem:lowerdiag}
Let $Y = X_{\lfloor d/2 \rfloor}$.  Each variable $x_{ij}$ appears exactly once on the union of the diagonal and the first and second superdiagonals of $Y$.  Let $(k,l)$ denote the indices of this unique occurrence in $Y$.  Then the matrix $Y^\searrow_{k+1,l+1}$ contains every variable in $\kk[x]$ with a lexicographically later index that is incomparable to $x_{ij}$ in $P_d$. 
\end{lemma}

\begin{proof}
That there is a unique occurrence of every variable follows from the fact that each of the matrices $A_i$ intersect the set of three super diagonals exactly once, and for each such $A_i$, only its superdiagonal or subdiagonal is intercepted (but not both).  To show the second claim, it suffices, by symmetry, to show this for the variables $x_{00}$, $x_{01}$, and $x_{02}$.  
The variables that are lexicographically later than $x_{00}$ and incomparable to it are $x_{02}$, $x_{11}$, $x_{12}$ and $x_{ij}$ for $i > 1$.  These variables all appear in 
$$Y^\searrow_{22} =\left(
\begin{array}{c| cc c}
x_{02} &  x_{11} & x_{12} & \cdots \\
\hline
x_{11} & x_{20} & x_{21} & \cdots \\
x_{12} & x_{21} & x_{22} & \cdots \\
\vdots & \vdots & \vdots & \ddots
\end{array} \right)$$
since the lower right block of this matrix is just the block Hankel matrix
$$\begin{pmatrix}
A_2 & A_3 & \cdots & A_{\lceil d/2 \rceil} \\
A_3 & A_4 & \cdots & A_{\lceil d/2 \rceil + 1} \\
\vdots & \vdots & \ddots & \vdots \\
A_{\lfloor d/2 \rfloor} & A_{\lfloor d/2 \rfloor+1} & \cdots & A_d 
\end{pmatrix}. $$
Similarly, the variables that are lexicographically later than $x_{01}$ are $x_{11}$, $x_{12}$ and all variables $x_{ij}$ with $i >1$.  But these all appear in $Y^\searrow_{23}$ which is obtained from deleting the first column of $Y^\searrow_{22}$. A similar argument shows the result for $x_{02}$.
\end{proof}

\begin{thm}
Any diagonal term order $\prec$ is delightful for the ideal $I_d$ for $2 \times 2$ minors of the block Hankel matrix $X_{\lfloor d/2 \rfloor}$.  The ideal $I_d$ is the prime ideal defining $\nu_d(\pp^1) \times \nu_2(\pp^1)$.  The $(r+1) \times (r+1)$ minors of $X_{\lfloor d/2 \rfloor}$ form a Gr\"obner basis for $I_d^{\{r\}}$ with respect to $\prec$.
\end{thm}

\begin{proof}
We must show that every $r+1$ element antichain of $P_d$ is the leading term of an $(r+1) \times (r+1)$ minor of $Y = X_{\lfloor d/2 \rfloor}$.  We begin by placing the elements of the antichain into ascending lexicographic order $ \{(i_0, j_0), \ldots, (i_{r}, j_{r}) \}$.  Thus, it suffices to show that there is a sequence of indices $(k_0, l_0), \ldots (k_r, l_r)$ such that $y_{k_t, l_t} = x_{i_t, j_t}$ for all $t$ and  $k_t< k_{t+1}$ and $l_t < l_{t+1}$ for all $t$.  We proceed by induction on $r$.

First of all, we can reduce to the  case where $i_0 = 0$. To see this, 
let $(k, l)$ be the unique index of $Y$ such that $y_{k, l} = x_{i_0, 0}$ and $y_{k, l}$ is on the main diagonal or the first or second superdiagonal.  The matrix $Y^{\searrow}_{k,l}$ has the form
$$Z=
\begin{pmatrix}
A_{i_0} &  A_{i_0 +1} & \cdots & A_{i_0 + \lceil d/2 \rceil + \sigma(i_0)} \\
A_{i_0 + 1} & A_{i_0 + 2} & \cdots & A_{i_0 + \lceil d/2 \rceil + 1 + \sigma(i_0)} \\
\vdots & \vdots & \ddots & \vdots \\
A_{i_0 + \lfloor d/2 \rfloor}  & A_{i_0 + \lfloor d/2 \rfloor+1 } & \cdots & A_d 
\end{pmatrix}
$$
where $\sigma(i_0) =  -1$ if $i_0$ is odd and $0$ otherwise.  If $i_0$ is even or if $i_0$ and $d$ are both odd, the matrix $Z$ has the form of $X_{\lfloor d_i/2 \rfloor}$ for some $d_i$.  If $i_0$ is odd and $d$ is even, $Z$ has the form $X_{\lceil d_i/2 \rceil} = X_{\lfloor d_i/2 \rfloor}^T$ for some $d_i$.  In any case, we may suppose that $i_0 = 0$, since all variables of interest appear inside a matrix of form $X_{\lfloor d_i/2 \rfloor}$, by Lemma \ref{lem:lowerdiag}.

Now if $i_1 = 0$ as well, we must have $(i_0, j_0) = (0,0)$ and $(i_1, j_1) =(0,2)$ and we set $(k_0,l_0) = (1,1)$ and $(k_1, l_1) = (2,2)$.  All remaining variables in the antichain lie in the matrix $Y^\searrow_{33}$ by Lemma \ref{lem:lowerdiag} which is of the form $X_{\lfloor (d-1)/2 \rfloor}$.  By induction, the rest of the antichain is a main diagonal of a minor which lies entirely within $Y^\searrow_{33}$. Thus the resulting sequence $(k_2, l_2), \ldots (k_r, l_r)$ satisfies $k_t > 2$ and $l_t > 2$.  Thus the concatenated sequence $(k_0, l_0),(k_1, l_1), \ldots (k_r, l_r)$ is a main diagonal sequence.  

If $(i_0, j_0) = (0,2)$  or if $i_1 > 1$ the argument is the same as the preceding paragraph.  The only remaining possibility is that the the sequence begins with one of the strings
$$(0,0), (1,1), (2,1), \ldots, (u-1),  (u,1) $$
$$(0,1),  (1,1), (2,1), \ldots,  (u-1), (u,1)  $$
$$(0,0), (1,1), (2,1), \ldots, (u-1),  (u,2) $$
$$(0,1),  (1,1), (2,1), \ldots,  (u-1), (u,2)  $$
and such that $i_{u+1} > u+1$.  In any of these cases, the beginning of the string is clearly a diagonal sequence by reading the unique elements on the main diagonal and the first and second superdiagonals of $Y$.  The condition that $i_{u+1} > u+1$ guarantees that all remaining variables in the sequence lie in lower right submatrix
$$
\begin{pmatrix}
A_{i_u+2} &  A_{i_u +3} & \cdots & A_{i_u + 2 \lceil d/2 \rceil + \sigma(u)} \\
A_{i_u + 3} & A_{i_u + 4} & \cdots & A_{i_u + 3  + \lceil d/2 \rceil +  \sigma(u)} \\
\vdots & \vdots & \ddots & \vdots \\
A_{i_u + 2 + \lfloor d/2 \rfloor}  & A_{i_u + 3 + \lfloor d/2 \rfloor } & \cdots & A_d 
\end{pmatrix}
$$
which as already shown, is of the form either $X_{\lfloor d_i/2 \rfloor}$ or $X_{\lfloor d_i/2 \rfloor}^T$.  In either case, by induction, the remaining part of the antichain is part of a diagonal sequence, the union with the diagonal sequence $(k_0, l_0) , \ldots, (k_u,l_u)$ will necessarily be a diagonal sequence.  This completes the proof.
\end{proof}

\begin{cor}\label{cor:blockhankel}
The ideal $I_d$ of $2 \times 2$ minors of the generic block Hankel matrix $X_{\lfloor d/2 \rfloor}$ is differentially perfect.  In particular, the set 
$$\left\{ \prod_{i =1}^l f^{i} \, \, | \, \,  f_{i}  \in \mathcal{M}_{r+1}(X_{\lfloor d/2 \rfloor}), \, \, \sum_{i =1}^l (\deg f_i - r)  = s  \right\}$$ 
 forms a Gr\"obner basis for the symbolic power $(I_d^{\{r\}})^{(s)}$ with respect to any diagonal term order.
 \end{cor}

\begin{proof}
The term order $\prec$ is delightful and the initial ideal is the antichain ideal of a poset.  The result follows by Theorem \ref{thm:gb}, Corollary \ref{cor:delightperfect}, and the fact that the incomparability graph of a poset is  a Greene graph.
\end{proof}


\subsection{The Scroll $\nu_d(\pp^1) \times \pp^k  $} \label{sec:scroll}

The Segre-Veronese varieties $ \nu_d(\pp^1) \times \pp^k $ are examples of scrolls and the delightfulness of a diagonal term order for the associated ideal of $2 \times 2$ minors was studied in Section 5 of \cite{Sturmfels2006}.  We wish to extend the construction described there to symbolic powers.  While the basic determinantal setup shares many features with the two preceding examples, one special feature here is that the initial ideal of the scroll is not the antichain ideal of any poset.   In all other cases where a straightening law and KRS correspondence is used, the corresponding initial ideal is an antichain ideal.  Thus, it is not clear that these standard techniques will work in this situation.

Let $\kk[x] := \kk[x_{ij} \, \, | \, \, i \in \{0,1,\ldots, d\}, j \in \{0,1, \ldots, k\} ]$ be the polynomial ring in $(d+1)(k+1)$ indeterminates.  For each $q \geq 0$ and $r\geq 1$  and let $X_{q,r}$ be the $(r+1) \times (k+1)$ generic matrix
$$X_{q,r} = 
\begin{pmatrix}
x_{q,0} &  x_{q,1} & \cdots & x_{q,k} \\
x_{q+1,0} & x_{q+1,1} & \cdots & x_{q+1,k} \\
\vdots & \vdots & \ddots & \vdots \\
x_{q+r,0} & x_{q+r,1} & \cdots & x_{q+r,k}
\end{pmatrix}
$$
and let $X_r$ be the concatenation of the $X_{q,r}$
$$X_r = \begin{pmatrix}
X_{0,r} & X_{1,r} & \cdots & X_{d-r,r} \end{pmatrix}. $$

Let $I_{d,k}$ be the ideal generated by the $2 \times 2$ minors of $X_1$  Let $\prec$ be any term order that selects the main diagonal of every minor of $X_r$ as the leading term.  One such term order is the lexicographic term order with $x_{i_1j_1}  \prec x_{i_2j_2}$ if $i_1 > i_2$ or if $i_1 = i_2$ and $j_1 > j_2$.

\begin{thm}{\rm \cite[Thm 5.9]{Sturmfels2006}} \label{thm:scrolls}
The $2 \times 2$ minors of $X_1$ are a Gr\"obner basis for $I_{d,k}$ with respect to the diagonal term order $\prec$.  This term order is delightful for $I_{d,k}$ and the $(r+1) \times (r+1)$ minors of $X_r$ form a Gr\"obner basis for $I_{d,k}$.
\end{thm}

To extend Theorem \ref{thm:scrolls} from secant ideals to symbolic powers, we must show that the quadratic initial ideal ${\rm in}_\prec(I_{d,k})$ is the edge ideal of a Greene graph. According to the proof of Theorem \label{thm:delscrolls} in \cite{Sturmfels2006}, this graph has vertices the pairs $(i,j)$ with $(i_1,j_1)$ connected to $(i_2,j_2)$ if $i_2 > i_1 +1$ or if $i_2 = i_1 + 1$ and $j_2 > j_1$.  This graph is not the incomparability graph of a poset, as was the case in all the preceding examples.  However, it turns out that it is the \emph{comparability graph} of a poset.  Indeed, define the poset $P_{d,k}$ on pairs $(i,j)$ subject to the relations $(i_1,j_1) \prec (i_2,j_2)$ if $i_2 > i_1 +1$ or if $i_2 = i_1 +1$ and $j_2 > j_1$.  This relation is clearly transitive, and hence defines a partial order.  The comparable pairs in $P_{d,k}$ correspond to the initial terms in the quadratic Gr\"obner basis for $I_{d,k}$.  Since the comparability graphs of posets are Greene graphs, we deduce:

\begin{cor}\label{cor:scroll}
The ideal $I_{d,k}$ of $2 \times 2$ minors of the matrix $X_1$ is differentially perfect.  In particular, the set 
$$\left\{ \prod_{i =1}^l f_{j_i}^{r_i} \, \, | \, \, f^{r_i}_{j_i} \in \mathcal{M}_{r_i + 1} (X_{r_i}), \, \, r_i \geq r,   \, \,  \sum_{i =1}^l (\deg f^{r_i}_{j_i} - r) = s \right\}$$ 
forms a Gr\"obner basis for the symbolic power $(I_{d,k}^{\{r\}})^{(s)}$ with respect to any diagonal term order.
\end{cor}


\end{document}